\documentclass[12pt]{article}
\everymath{\displaystyle}
\usepackage{graphicx}
\usepackage{pst-all}
\usepackage{color}
\usepackage{amssymb}
\usepackage{amsmath}

\newtheorem{prethm}{{\bf Theorem}}

\newenvironment{thm}{\begin{prethm}{\hspace{-0.5
               em}{\bf.}}}{\end{prethm}}
\newtheorem{prelemma}{{\bf Lemma}}

\newtheorem{preex}{{\bf Example}}

\newtheorem{preprop}{{\bf Proposition}}

\newenvironment{prop}{\begin{preprop}{\hspace{-0.5em}{\bf.}}}{\end{preprop}}
\newtheorem{precor}{{\bf Corollary}}

\newenvironment{cor}{\begin{precor}{\hspace{-1
               em}{\bf.}}}{\end{precor}}
\newtheorem{preremark}{{\bf Remark}}

\newtheorem{preprob}{{\bf Problem}}

\newtheorem{predefin}{{\bf Definition}}

\newtheorem{preconj}{{\bf Conjecture}}

\newtheorem{preprobb}{{\bf Problem}}

\newtheorem{prelem}{{\bf Theorem}}

\newenvironment{proof}{{\bf Proof.}\rm }{\hfill{$\Box$}}

\newtheorem{presolution}{{\bf Solution.}}

\def\newpic#1{}

\title{\Large\bf Bounds for the $b$-chromatic number of some families of graphs}
\author{
{\large {\bf Mekkia Kouider}$^{\bf 1}$ and \large {\bf Manouchehr
Zaker}$^{\bf 2}$ \footnote{while visiting university of Orsay}}
\vspace{3mm}\\
    $^{\bf 1}$ L.R.I, Unite Mixte de Recherche 8623\\
    Universit\'e Paris-Sud, Orsay, France\\
    km@lri.fr.\\
    {}\\
    $^{\bf 2}$ Institute for Advanced Studies in Basic Sciences\\
    45195-1159, Zanjan, Iran\\
    mzaker@iasbs.ac.ir.}

\date{December 2004}
\begin{document}
\maketitle

\begin{abstract}
In this paper we obtain some upper bounds for $b$-chromatic number
of $K_{1,t}$ -free graphs, graphs with given minimum clique
partition and bipartite graphs. These bounds are in terms of
either clique number or chromatic number of graphs or biclique
number for bipartite graphs. We show that all the bounds are
tight.
\end{abstract}

\section{Introduction and related results}
All graphs in this paper are finite, simple and undirected graphs.
By clique number of a graph $G$ we mean the largest order of a
complete subgraph in $G$ and denote it by $\omega(G)$. Also
$\alpha(G)$ stands for the largest number of independent vertices
in $G$. For other notations which are not defined here we refer
the reader to \cite{Bondy}.

An {\it antimatching} of a graph $G$ is a matching of its complement.

A {\it proper coloring} of $G$ is a coloring of the vertices such
that any two adjacent vertices have different colors. Given a
proper coloring of $G$, a {\it $t$-dominating set} $T = \{ x_1
,...,x_t \}$ is a set of vertices such that $T$ is colored by $t$
colors and each $x_i$ is adjacent to $t-1$ vertices of different
colors. In that case, and if $G$ is colored by exactly $t$ colors,
we say we have a {\it $t$-dominating coloring} (or {\it
$b$-coloring} with $t$ colors). We denote by $\varphi(G)$ the
maximum number $t$ for which there exists a $t$-dominating set in
a coloring of $V(G)$ by $t$ colors. This parameter has been
defined by Irving and Manlove \cite{IM}, and is called the {\it
$b$-chromatic number} of $G$. In a $b$-coloring of a graph $G$
with $b$ colors, any vertex $v$ which has at least $b-1$ neighbors
with different colors is called a {\it representative}. We note
that in any $b$-coloring of $G$ with $b$ colors there should be at
least $b$ representatives with $b$ different colors.

It is known that $\chi(G) \leq \varphi(G) \leq \Delta~+1$.
Let $G$ be a graph with decreasing degree sequence $d(x_1) \geq
d(x_2) \geq \ldots \geq d(x_n)$ and let $m(G)= max~ \{i : d(x_i)
\geq (i-1) \}$. In \cite{IM}, the authors proved that for any
graph $G$, $\varphi(G) \leq m(G)$ and they show that for tree $T$
the inequality $m(T)-1 \leq \varphi(G) \leq m(T)$ is satisfied.

Also in \cite{IM} it is shown that determining $\varphi$ is
NP-hard for general graphs, but polynomial for trees.

Some authors have obtained upper or lower bounds for $\varphi(G)$
when $G$ belongs to some special families of graphs. In \cite{K},
$b$-chromatic number of graphs with girths five and six has been
studied. Let $G$ be a graph of girth at least $5$, of minimum
degree $\delta$ and of diameter $D$, it is shown in \cite{K} that
$\varphi(G)> min \{\delta, D/6\}$ and that if $G$ is $d$-regular,
of girth at least six, then $\varphi(G) = d+1$. In this last case
the construction of a $b$-dominating coloring is done in a
polynomial time.

Kratochvil et al. in \cite{KTV} showed that for a $d$-regular
graph $G$ with at least $d^4$ vertices, $\varphi(G)=d+1$.

In \cite{KM}, Kouider and Mah\'eo discuss on the $b$-chromatic
number of the cartesian product $G \Box H$ of two graphs $G$ and
$H$. They prove that $\varphi (G \Box H)\geq \varphi (G)+ \varphi
(H)-1$ when $G$ (resp. $H$) admits $\varphi(G)$ (resp.
$\varphi(H)$) dominating set which is stable set.

We also recall the following result of Klein and Kouider
\cite{KK}. Let $\mathcal D$ be $K_4 \setminus e$. Let $G$ be a
$P_4$-free graph, then $\varphi(G)= \omega(G)$, for any induced
subgraph of $G$ if and only if $G$ is $2 \mathcal D$-free and
$3P_3$-free.

The aim of this paper is to obtain an upper bound for
$b$-chromatic number of a graph $G$ when $G$ is restricted to be
in special families of graphs. In section $2$ we consider
$K_{1,t}$ -free graphs. In section $3$ we give an upper bound in
terms of clique number and minimum clique partition of a graph.
Finally in section $4$ bipartite graphs will be considered. We
also show that all the bounds obtained in this paper are tight.

\section{$K_{1,t}$ -free graphs}

In this section we give an upper bound for the $b$-chromatic
number of $K_{1,t}$ -free graphs, when $t\geq 3$. If $t=2$ then
the graph should be a complete graph for which the $b$-chromatic
number is the same as chromatic number.

\begin{thm}\label{free}
Let $G$ be a $K_{1,t}$ -free graph where $t\geq 3$, then
$\varphi(G) \leq (t-1)(\chi(G)-1)+1$.
\end{thm}

\begin{proof}
Suppose $\varphi(G)=b$. Let $C$ be a color class in a $b$-coloring
of $G$ with $b$ colors, and let $x$ be any representative of the
class $C$. Among the neighbors of the vertex $x$ there exist a set
say $S$ of $b-1$ vertices with distinct colors. Let $H$ be the
subgraph induced by $S$. By the hypothesis on the graph $G$ we
have $\alpha(H) \leq t-1$ and also $\chi(H)\leq \chi(G) - 1$. So
$$b-1=|V(H)|\leq \alpha(H).\chi(H)\leq (t-1)(\chi(G)-1).$$ Therefore
$b\leq (t-1)(\chi(G)-1)+1$.
\end{proof}

In the following we show that the bound of the theorem can be
achieved for each $t$.

\begin{prop}
For any integer $t\geq 3$ and $k$, there exists a $K_{1,t}$ -free
graph $G$ such that $\chi(G)=k$ and $\varphi(G) = (t-1)(k-1)+1$.
\end{prop}

\begin{proof}
Suppose the graph $H$ is defined as a vertex $v$ such that its
neighbors form $t-1$ mutually disjoint cliques with $k-1$ vertices.
Now we take $(t-1)(k-1)+1$ disjoint copies of $H$ and connect them
sequentially by exactly one edge between any two consecutive copies.
These edges can be incident to any vertex other than $v$ and its
copies in other copies of $H$. We denote the resulting graph by $G$.
It is easily seen that $G$ satisfies the conditions of theorem.
\end{proof}

We have now the following immediate corollary of theorem \ref{free}.

\begin{cor}
If $G$ is a claw-free graph, then $\varphi(G)\leq 2\chi(G)-1$.
\end{cor}

In \cite{CS} the important fact $\chi(G)\leq 2\omega(G)$ is proved
for a claw-free graph $G$ satisfying $\alpha(G)\geq 3$, therefore
using this result we obtain $\varphi(G)\leq 4\omega(G) - 1$.

\section{$b$-coloring and minimum clique partition}

In this section we give a bound for the $b$-chromatic number of a
graph $G$ in terms of its minimum clique partition. A clique
partition for a graph $G$ is any partition of $V(G)$ into subsets
say $C_1, C_2, \ldots, C_k$ in such a way that the subgraph of $G$
induced by $C_i$ is a clique, for each $i$. We denote by
$\theta(G)$ the minimum number of subsets in a clique partition of
the graph $G$. We note that for any graph $G$, $\chi(\overline{G})
= \theta(G)$; also, if $\theta(G)=k$ then $G$ is the complement of
a $k$-partite graph. Therefore the following result applies for
all graphs.

\begin{thm}\label{clique-partition}
Let $G$ be a graph with clique partition number
$\theta(G)=k$ and clique number $\omega$,
then $\varphi(G) \leq \frac{k^2\omega}{2k-1}$.
\end{thm}

\begin{proof}
If $k=1$ then $G$ is complete and equality holds in the inequality
of theorem. We suppose now $k\geq 2$. As $\theta(G)=k$, therefore
$\alpha(G)\leq k$. Let us consider a $b$-coloring of $G$ with
$\varphi(G)=b$ colors. Let $i_j$ be the number of color classes
with exactly $j$ elements. As  $\alpha(G)\leq k$, we know that
$i_j=0$ for $j\geq k+1$. So we have
$$b=\sum_{j=1}^{k}i_j.$$

By hypothesis, there exists a partition of $V(G)$ into $k$
complete subgraphs, therefore if $n$ is the order of $G$,
$$n = \sum_{j=1}^{k}j.i_j =b+ \sum_{j=2}^{k}(j-1)i_j\leq
k\omega.~~~~{\bf (1)}$$

Suppose first that $i_1=0$. Then any color class in the
$b$-coloring of $G$ with $b$ colors contains at least two
vertices. This shows that $b\leq n/2$ and so $b\leq k\omega/2$.
Finally $b \leq \frac{k^2}{2k-1}\omega$, because $\frac{k}{2}\leq
\frac{k^2}{2k-1}$.

Suppose now $i_1\geq 1$ and let $C_i=\{x_i \}$ for $i = 1, \ldots,
i_1$. Then any representative of any color $j$ is adjacent to any
$x_i$, where $i,j \leq i_1$ and $i \neq j$.

It follows that $\{x_1, \ldots, x_{i_1}\}$ induces a complete
subgraph of $G$. On the other hand, by the fact that there exists
a partition of $V(G)$ into $k$ cliques and the pigeonhole
principle, at least $\frac{\sum_{j=2}^{k}i_j}{k}$ of
representative vertices form a complete graph. We know from above
that any representative of any color $j$ is adjacent to any $x_i,
i \neq j,i \leq i_1 $, consequently there is a complete subgraph
of at least $i_1 ~+ \frac{\sum_{j=2}^{k}i_j}{k}$ vertices. We get
the following inequality
$$i_1 ~+ \frac{\sum_{j=2}^{k}i_j}{k} \leq \omega$$
in other words, $$ki_1 ~+ \sum_{j=2}^{k}i_j \leq k\omega.
\hspace{1cm} {\bf (2)}$$

Now we have
$$(2k-1)b=\sum_{j=1}^{k}(2k-1)i_j=(k-1)(\sum_{j=1}^{k}ji_j) + ki_1 + i_2
- \sum_{j=3}^{k}((k-1)j-2k+1)i_j~~for~ k\geq 3,$$ or
$$(2k-1)b=(k-1)(\sum_{j=1}^k ji_j)+ki_1+i_2~~ for~ k=2.$$

So we have $$(2k-1)b \leq (k-1)n + ki_1 + i_2$$ and by inequality
(1), $$(2k-1)b \leq (k-1)k\omega + ki_1 + i_2 \leq k^2 \omega - k(
\omega - i_1 - \frac{i_2}{k}).$$ By inequality (2),
$$(2k-1)b \leq k^2\omega.$$ The theorem is proved.
\end{proof}

\begin{preprop}
For any positive integers $k\geq 2$ and $\omega$ divisible by
$2k-1$, there exists a
graph $G$ with $\theta(G)=k$ and with clique number $\omega$, such
that $\varphi(G) = \frac{k^2\omega}{2k-1}$.
\end{preprop}

\begin{proof}
In order to construct our graph we first consider three sets of
mutually disjoint cliques $\{A_1,\ldots,A_k\}$,
$\{B_1,\ldots,B_k\}$ and $\{C_1,\ldots,C_k\}$ where $|A_i| =
\frac{\omega}{2k-1}$, $|B_i| = |C_i| = \frac{(k-1)\omega}{2k-1}$,
for each $i=1,\ldots,k$. We put an edge between any two vertices
$u$ and $v$ in $A_i$ and $A_j$ for each $i$ and $j$, therefore
$\bigcup_i A_i$ forms a clique with $\frac{k\omega}{2k-1}$
vertices. Then we join any vertex in $A_i$ to any vertex in $B_j$
for each $i$ and $j$, and also we join the vertices of $A_i$ to
all the vertices of $C_i$, for each $i$. We don't have any edge
between any two vertices of $B_i$ and $B_j$ when $i\neq j$ and the
same holds for $C_i$'s. Finally we put an edge between any two
vertices $v\in B_i$ and $u\in C_j$ if $i\neq j$.

We color the vertices in $\bigcup_i A_i$ with $1,2, \ldots,
\frac{k\omega}{2k-1}$ and the vertices of $\bigcup_i B_i$ with
distinct colors $\frac{k\omega}{2k-1}+1, \ldots,
\frac{k^2\omega}{2k-1}$. The colors in $C_i$ will be the same as
$B_i$ for each $i$. All the vertices of $A=\bigcup_i A_i$ are
representatives and the same holds for $B=\bigcup_i B_i$.

Now it is enough to show that the constructed graph $G$ has the
clique number $\omega$. We first observe that if we identify each
of cliques $A_i$'s, $B_i$'s and $C_i$'s with single vertices
$a_i$'s, $b_i$'s and $c_i$'s, respectively, then we may define a
graph $H$ with $3k$ vertices with vertex set
$\{a_1,\ldots,a_k,b_1,\dots,b_k,c_1,\ldots,c_k\}$ where there is
an edge between two vertices $u$ and $v$ if and only if their
corresponding cliques are jointed in the graph $G$. Therefore to
find the maximum number of vertices in a clique of the graph $G$,
it is enough to check all cliques in $H$. Let us first set $A=
\{a_1,\ldots,a_k\}$, $B=\{b_1,\dots,b_k\}$ and $C=
\{c_1,\ldots,c_k\}$. Let $K$ be a clique in $H$. There are two
possibilities:

{\bf 1}. There is no vertex from $C$ in $K$. In this case $K$ may
contain all vertices in $A$ and at most one from $B$, i.e. with at
most $k+1$ vertices. This clique results in a clique in $G$ with
$\frac{k\omega}{2k-1} + \frac{(k-1)\omega}{2k-1}=\omega$ vertices.

{\bf 2}. There is one vertex from $C$ in $K$. In this case $K$
contains only one vertex from $C$ and at most one vertex from $A$
and one from $B$. And this may happen when we consider for example
$a_1$ and its neighbor in $C$ and a suitable vertex in $B$. This
clique of order three results in a clique in $G$ with
$\frac{(k-1)\omega}{2k-1}+ \frac{\omega}{2k-1}=\omega$ vertices.
\end{proof}

The following result is an immediate corollary of theorem
\ref{clique-partition}.

\begin{cor}
For any graph $G$, with clique-number $\omega(G)$,
$$\varphi(G) \leq \frac{\chi^2(\overline{G})}{2~
\chi(\overline{G})-1}~\omega(G).$$
\end{cor}

In the case that $G$ is the complement of a bipartite graph we
have more knowledge on its $b$-colorings. We first introduce some
special graphs which play an important role in $b$-colorings of
the complement of bipartite graphs. Before we begin let us mention
that when we say there is an anti-matching between two subsets $X$
and $Y$ in a graph $G$, it means that there exists a matching
between $X$ and $Y$ in the complement of $G$.

Let $G$ be the complement of a bipartite graph with a bipartition
$(X,Y)$ in such a way that there are partitions of $X$ and $Y$
into three subsets as $X=A_1\cup B_1\cup C_1$ and $Y=A_2\cup
B_2\cup C_2$ such that the following properties hold:

\noindent {\bf 1.} Any vertex in $A_1$ is adjacent to any vertex
in $A_2\cup B_2$, hence the subgraph induced by $A_1\cup A_2 \cup
B_2$ in $G$ is a clique. Also any vertex in $A_2$ is adjacent to
any vertex in $C_1$.

\noindent {\bf 2.} $|B_1|=|B_2|$ and there is a perfect anti-matching
between $B_1$ and $B_2$.

\noindent {\bf 3.} $|C_1|=|C_2|$ and there is a perfect anti-matching
between $C_1$ and $C_2$.

In this case by letting $b=|A_1\cup A_2| + |B_1| +
|C_1|=|X|+|A_2|$, we say $G$ belongs to the family
$\mathcal{A}_b$. In fact $\mathcal{A}_b$ consists of all the
complement of bipartite graphs $G$ which admits the
above-mentioned properties. Let us remark that $\varphi(G) \geq b$
for any graph $G$ belonging to $\mathcal{A}_b$. In fact, we color
$X \cup A_2$ with different colors; using the antimatchings, we
give to $B_2$ the same colors as $B_1$, and to $C_2$ the same
colors as $C_1$.

\begin{thm}\label{cobip}
Let $G$ be the complement of a bipartite graph, then
$\varphi(G)\leq \frac{4\omega}{3}$. Furthermore, there is a
$b$-coloring for $G$ with $b$ colors if and only if $G$ is in
$\mathcal{A}_b$.
\end{thm}

\begin{proof}
The inequality $\varphi(G)\leq \frac{4\omega}{3}$ follows from
theorem \ref{clique-partition} where we put $k=2$.

If $G$ is in $\mathcal{A}_b$ then by the comment before theorem
\ref{cobip} there is a $b$-coloring for $G$ with $b$ colors.

Suppose now we have a $b$-coloring for $G=(X,Y;E)$ with $b$ colors
$\{1,2,\ldots,b\}$. Let the color classes be $U_1$, $U_2$, \ldots,
$U_b$ and without loss of generality we may suppose that $|U_i|=1$
for $i=1,\ldots, t$. Therefore $|U_i| = 2$ for $i
>t$. Set $A_1 = X \cap \cup_{i=1}^t U_i$ and $A_2 = Y \cap
\cup_{i=1}^t U_i$.

Let $u_i, i=t+1, ...,s$ be the representatives contained in $X$
and they form a set $B_1$; let $u_i, i=s+1,...,b$ be the remaining
representatives, these are by definition in $Y \setminus A_2$ and
they form a set $C_2$. As each color class for $i \geq t+1$, has
exactly $2$ elements, there exists a set $B_2$ in $Y$ with $|B_2|
= |B_1|$ with the same colors as $B_1$. Similarly there exists a
set $C_1$ in $Y$ with $|C_2| = |C_1|$ with the same colors as
$C_2$.

There are perfect anti-matchings, one between $B_1$ and $B_2$ and
another between $C_1$ and $C_2$. By the property of being
representative for each element of $B_1 \cup C_2$, and by the
unicity of the elements colored by the colors of $A_1 \cup A_2$,
$A_1 \cup C_2$ is a clique, $A_2 \cup B_1$ is also a clique.
Considering now the partitions $X=A_1\cup B_1\cup C_1$ and
$Y=A_2\cup B_2\cup C_2$ we conclude that $G$ belongs to
$\mathcal{A}_b$.
\end{proof}

We get easily the following consequence.

\begin{cor}
Let $G$ be the complement of a bipartite graph. Then

$$\varphi(G)= b$$ if and only if~
$\max \{k : G \in \mathcal{A}_k  \}~=b$.
\end{cor}

\vspace{5mm}

Let us remark that for the larger class of  graphs $G$ with
$\alpha(G)=2$, there is no linear bound for $b$-chromatic number
(even for chromatic number) in terms of $\omega(G)$ because, as
pointed out in \cite{CS}, or each $k$ there is a graph $G$ with
$\alpha(G)=2$ such that $\chi(G)\geq k/2$ and $\omega(G)=o(k)$.

\section{Bipartite graphs}

In this chapter we suppose $G$ is a bipartite graph. In the
following by the {\it biclique number} of $G$ we mean the minimum
number of mutually disjoint complete bipartite graphs which cover
the vertices of $G$. Any subgraph of $G$ which is complete
bipartite graph is called a {\it biclique} of $G$.

\begin{thm}\label{bipart}
Let $G$ be a bipartite graph with bipartition $(X,Y)$, on $n$
vertices and biclique number $t$. Then
$$\varphi(G)\leq  \lfloor \frac{n-t+4}{2} \rfloor.$$
\end{thm}

\begin{proof}
We first prove the theorem for graphs $G=(X,Y)$ which admits a
$b$-coloring with $b=\varphi(G)$ colors where there is at least
one representative in $X$ and also one in $Y$. Let these
representatives be $v\in X$ and $u\in Y$. Then $v$ has at least
$b-1$ neighbors in $Y$ and also $u$ has at least $b-1$ neighbors
in $X$. These give us two bicliques with cardinality at least
$2b-2$ and at most $2b$. As $t$ is the biclique number of $G$
there should be at least $t-2$ vertices in $G$. Therefore $n\geq
2b-2+t-2$ and so $b\leq \frac{n-t+4}{2}$.

Now we may suppose that in a $b$-coloring of $G$, all the
representatives are in a same part say $X$. Let $i_j$ be the number
of color classes in the $b$-coloring with exactly $j$ colors in part
$Y$. There are two possibilities.

Suppose first that $i_1\geq 1$. Let $w$ be the vertex in any color
class with cardinality one in the part $Y$. Then $w$ belongs to
$Y$ and has $b-1$ neighbors which are representatives of different
colors. So $w$ is representative. This is a contradiction with the
hypothesis on $X$.

Now let $p$ be the minimum number with $i_p\neq 0$. So $p\geq 2$.
We have $n\geq b+ bp = b(p+1)$. We may suppose at this stage that
all vertices in $X$ are representatives and of different colors,
and, also any vertex $y$ in $Y$ is adjacent to, at least, some
representative and is the unique vertex of color c(y) of this
representative. Otherwise, if we delete those vertices in $X$
which are not representatives and also vertices in $Y$ without the
previous property, we prove the inequality of the theorem, for the
resulting graph $G'$. Let $n-l$ be its order. We have
$$\varphi(G')\leq \frac{n-l-t'+4}{2}.$$
As the inequality $t\leq t'+l$  holds, it can be seen easily that
we get $$ \varphi(G')\leq \frac{n-t+4}{2}.$$ We also have, by
construction of $G'$, $\varphi(G)\leq \varphi(G')$.

So it is enough to prove the theorem for the case where all the
vertices in $X$ are representative and any vertex in $Y$ is
adjacent to some representative. By these hypothesis, and as the
coloring is proper, we have $t\leq b$. Finally since $n\geq
b(p+1)$ and $p\geq 2$ then $2b \leq  n - (p-1)b \leq n-b \leq
n-t$. Therefore $b\leq \frac{n- t}{2}$.
\end{proof}

\begin{prop}
For any integer $p\geq 3$, there is a bipartite graph $G$ with
$n=3p-4$ vertices and biclique number $t=p-1$ such that $b(G) = p =
\lfloor \frac{n-t+4}{2} \rfloor$.
\end{prop}

\begin{proof}
We first consider a complete bipartite graph $K_{p-1,p-1}$ minus a
matching with size $p-2$. We color one part say $X$ of this graph
with $1,3,4, \ldots, p$ and other part say $Y$ with $2,3,4,
\ldots, p$ so that vertices with colors $1$ and $2$ are adjacent.
Then we add $p-2$ extra vertices to the part $X$ and color all of
them with $2$. Now put a matching with size $p-2$ between these
extra vertices in $X$ and all the vertices in $Y$ except the one
colored by $2$. The resulting graph $G$ is a graph of order
$n=3p-4$ with a $b$-coloring with $p$ colors. In fact, $b(G)$ is
exactly equal to $p$ because $\Delta(G)=p-1$. By the precedent
theorem,
$$b(G) \leq \frac{n-t+4}{2}.$$

It is then enough to show that $t \geq p-1$ to get the reverse
inequality. Because there are $p-2$ vertices with degree one, at
least $p-2$ cliques are required to cover these vertices. We
observe that we need an extra clique to cover the vertex colored
by $2$ in $Y$. Now we get the equality $b(G) = \lfloor
\frac{n-t+4}{2} \rfloor$.
\end{proof}


\begin{thebibliography}{1}
\bibitem{Bondy}
J. A. Bondy and U. S. R. Murty, {\it Graph Theory with Applications},
American Elsevier Publishing Co., Inc., 1976.

\bibitem{CS}
M. Chundovsky and P. D. Seymour, {\it The structure of claw-free
graphs}, manuscript 2004.

\bibitem{IM}
I. W. Irving and D. F. Manlove, {\it The $b$-chromatic number of a
graph, Discrete Applied Math.}, 91 (1999), 127-141.

\bibitem{KK}
S. Klein and M. Kouider, {\it $b$-coloration and $P_4$-free graphs},
manuscript.

\bibitem{K}
M. Kouider, {\it b-chromatic number, subgraphs and degrees},
submitted.

\bibitem{KM}
M. Kouider and M. Mah\'eo, {\it Some bounds for the $b$-chromatic
number of a graph}, Discrete Math., 256 (2002), 267--277.

\bibitem{KTV}
J. Kratochvil, Zs. Tuza, and M. Voigt, {\it On the $b$-chromatic
number of graphs}, Lecture Notes in Computer Science, Springer,
Berlin, 2573 (2002), 310--320.

\end{thebibliography}
\end{document}